\documentclass{tran-l}

\usepackage{amsthm}
\usepackage{amsmath}
\usepackage{amssymb}
\usepackage{latexsym}
\usepackage{enumerate}
\usepackage{graphics} 
\usepackage{graphicx}

\newcommand{\grad}{{\rm{grad}\,}}
\newcommand{\Hess}{{\rm Hess}\,}

\newcommand{\hessian}{{{\rm Hess}}}
\newcommand{\gradien}{{\rm{grad}}}
\newcommand{\R}{\mathbb{R}}

\newtheorem{theorem}{Theorem}
\newtheorem{lemma}[theorem]{Lemma}

\newtheorem{corollary}[theorem]{Corollary}

\theoremstyle{definition}

\newtheorem{remark}[theorem]{Remark}

\numberwithin{equation}{section}

\newcommand{\rf}[1]{\mbox{(\ref{#1})}}
\newcommand{\rl}[1]{{~\ref{#1}}}

\def\beq{\begin{equation}}
\def\eeq{\end{equation}}

\def\qed{\ifhmode\unskip\nobreak\fi\ifmmode\ifinner\else
\hskip5 pt \fi\fi\hbox{\hskip5 pt \vrule width4 pt height6 pt
depth1.5 pt \hskip 1pt }}

\begin{document}

\title[The sectional curvature of  cylindrically  bounded submanifolds]{An estimate for the sectional curvature of  cylindrically  bounded submanifolds}

\author[L. J. Al\'\i as]{Luis J. Al\'\i as$^{1}$}
\address{Departamento de Matem\'{a}ticas, Universidad de Murcia, E-30100 Espinardo, Murcia, Spain}
\email{ljalias@um.es}
\author[G.P. Bessa]{G.Pacelli Bessa$^{2}$}
\address{Departamento de Matem\'{a}tica, Universidade Federal do Ceara-UFC, Bloco 914, Campus do Pici, 60455-760, Fortaleza, Ceara, Brazil}
\email{bessa@mat.ufc.br}

\author[J.F. Montenegro]{J. Fabio Montenegro$^{3}$}
\address{Departamento de Matem\'{a}tica, Universidade Federal do Ceara-UFC, Bloco 914, Campus do Pici, 60455-760, Fortaleza, Ceara, Brazil}
\email{fabio@mat.ufc.br}

\thanks{$^{1}$ This research is a result of the activity developed within the framework of the
Programme in Support of Excellence Groups of the Regi\'{o}n de Murcia, Spain, by Fundaci\'{o}n S\'{e}neca,
Regional Agency for Science and Technology (Regional Plan for Science and Technology 2007-2010).
Research partially supported by MEC project PCI2006-A7-0532,
MICINN project MTM2009-10418, and Fundaci\'{o}n S\'{e}neca project 04540/GERM/06,
Spain.\\
\indent $^{2}$ Partially supported by CNPq-CAPES (Brazil) and MEC project PCI2006-A7-0532 (Spain).\\
\indent $^{3}$ Partially supported by CNPq-CAPES (Brazil).}

\subjclass[2000]{53C42}

\keywords{Omori-Yau Maximum Principle, cylindrically bounded submanifolds, properly immersed submanifolds.}

\date{June 2010}

\begin{abstract}
We give sharp sectional curvature estimates for complete immersed  cylindrically bounded $m$-submanifolds
$\varphi\colon M^m \to N^{n-\ell}\times \mathbb{R}^{\ell}$, $n+\ell \leq 2m-1$ provided that either $\varphi $ is proper with the
norm of the second fundamental form with certain controlled growth or $M$ has  scalar curvature with strong quadratic decay. This latter
gives a non-trivial extension of the Jorge-Koutrofiotis Theorem \cite{jorge-koutrofiotis}. In the particular case of
hypersurfaces, that is, $m=n-1$, the growth rate of the norm of the second fundamental form is improved. Our results will be an
application of a generalized Omori-Yau Maximum Principle for the Hessian of a Riemannian manifold, in its newest elaboration given
by Pigola, Rigoli and Setti in \cite{pigola-rigoli-setti}.
\end{abstract}

\maketitle

\section{Introduction}\hspace{.5cm} Given  complete  Riemannian manifolds $M^m$ and  $N^n$ with dimension $m<n$, the isometric immersion problem asks whether there exists an isometric immersion
$\varphi \colon M\hookrightarrow  N$. When $N^{n}=\mathbb{R}^{n}$ is the Euclidean space, the isometric problem is answered by  the Nash Embedding Theorem that  says that there is  an isometric embedding $\varphi \colon M^{m}
\hookrightarrow \mathbb{R}^{n}$  provided  the codimension $n-m$ is sufficiently large, see \cite{nash}. For small codimension, meaning in this paper that $n-m\leq m-1 $,  the answer in general depends on the geometries of $M$ and $N$. For
instance, the Hilbert-Efimov Theorem  \cite{efimov}, \cite{hilbert} says that no complete surface $M$ with sectional curvature $K_{M}\leq -\delta^{2}<0$ can be isometrically immersed in $\mathbb{R}^{3}$ and a classical result by C.
Tompkins \cite{tompkins}, states that a compact, flat, $m$-dimensional Riemannian manifold can not be isometrically immersed in $\mathbb{R}^{2m-1}$. Tompkins' result was extended in  a series of papers, by Chern and Kuiper
\cite{chern-kuiper}, Moore \cite{moore}, O'Neill \cite{oneil},  Otsuki \cite{otsuki} and Stiel \cite{stiel},  whose results  can
be summarized in the following theorem (we recall that a Cartan-Hadamard manifold is a simply connected, complete, Riemannian manifold
with non-positive sectional curvatures).

\begin{theorem} \label{thmTompkins} Let  $\varphi\colon M^{m}\hookrightarrow N^{n}$, $n\leq 2m-1$, be an isometric immersion of a compact  Riemannian $m$-manifold $M$ into  a Cartan-Hadamard $n$-manifold $N$.
Then the sectional curvatures of $M$ and $N$ satisfy \begin{equation}\sup_{M} K_{M}> \inf_{N} K_{N}. \label{eqTompkins}\end{equation}
\end{theorem}

Theorem \ref{thmTompkins} was extended by  Jorge and Koutrofiotis in  \cite{jorge-koutrofiotis} to
bounded, complete submanifolds with scalar curvature bounded below and in the version presented by Pigola, Rigoli and Setti
in \cite[Theorem 1.15]{pigola-rigoli-setti}, with scalar curvature satisfying
\begin{equation}
s_{M}(x)\geq -B^2\varrho^{2}_{M}(x)\cdot\prod_{j=1}^{k}\Big(\log^{(j)}(\varrho_{M}(x))\Big)^2,\,\, \varrho_{M}(x)\gg 1,
\label{eqScalar}
\end{equation}
for some constant $B>0$ and some integer $k\geq 1$,
where $\varrho_{M}$ is the distance function on $M$ to a fixed point and $\log^{(j)}$ is the $j$-th iterate of the
logarithm.

\begin{theorem}[Jorge-Koutrofiotis, \cite{jorge-koutrofiotis}]\label{thmJorge-koutrofiotis}
Let $M^{m}$ and $N^{n}$  be complete Riemannian manifolds of dimensions $m$ and $n$, respectively, with $n\leq 2m-1$ and
let $\varphi:M\rightarrow N$ be an isometric immersion with $\varphi(M)\subset B_N(r)$, where $B_{N}(r)$
denotes a geodesic ball of $N$ centered at a point $p\in N$ and radius $r$.
Assume that the radial sectional curvature $K_N^{\mathrm{rad}}$ along the radial geodesics issuing from $p$
satisfies $K_N^{\mathrm{rad}}\leq b$ in $B_N(r)$ and $0<r<\min\{\mathrm{inj}_N(p),\pi/2\sqrt{b}\}$, where we replace
$\pi/2\sqrt{b}$ by $+\infty$ if $b\leq 0$.
If the scalar curvature of $M$ satisfies (\ref{eqScalar}), then
\begin{equation}\label{eqJorge-koutrofiotis}
\sup_{M}K_{M}\geq C_{b}^{2}(r)+\inf_{B_{N}(r)}K_{N},
\end{equation}
where
\[
C_b(t)=
\begin{cases}
\sqrt{b}\cot(\sqrt{b}\,t) & \mbox{\rm if $b>0$ and $0<t<\pi/2\sqrt{b}$}\\
1/t & \mbox{\rm if $b=0$ and $t>0$}\\
\sqrt{-b}\coth(\sqrt{-b}\,t) & \mbox{\rm if $b<0$ and $t>0$}.
\end{cases}
\]
\end{theorem}

\begin{remark} If $N^{n}=\mathbb{N}^{n}(b)$ is the simply connected space form  of constant sectional curvature $b$ and
$M=\partial B_{\mathbb{N}^{n}(b)}(r)\subset \mathbb{N}^{n}(b)$ is a geodesic sphere of radius $r$ then equality in
(\ref{eqJorge-koutrofiotis}) is achieved.
\end{remark}

The purpose of this paper is to extend Jorge-Koutrofiotis theorem to the case of complete cylindrically bounded
submanifolds of a Riemannian product $N^{n-\ell}\times\mathbb{R}^{\ell}$, where $N$ is a complete Riemannian manifold
of dimension $n-\ell$. In this context, an isometric immersion $\varphi:M^m\rightarrow N^{n-\ell}\times\mathbb{R}^{\ell}$
of a Riemannian manifold $M^m$ is said to be \textit{cylindrically bounded} if there exists $B_N(r)$, a geodesic ball of
$N$ centered at a point $p\in N$ with radius $r>0$, such that $\varphi(M)\subset B_N(r)\times\mathbb{R}^{\ell}$ (see Figure 1).
\begin{figure}[h]
\begin{center}
\includegraphics[width=3cm]{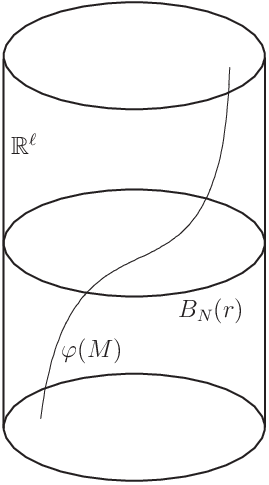}
\caption{A cylindrically bounded submanifold.}
\end{center}
\end{figure}

In a recent paper \cite{ABD}, the two first authors, jointly with Dajczer, derived an estimate for the mean curvature $H$ of complete
cylindrically bounded submanifolds into a product Riemannian  manifold $N^{n-\ell}\times\R^{\ell}$. Specifically, it was proved in \cite{ABD}
that if the immersion $\varphi:M^m\rightarrow N^{n-\ell}\times\mathbb{R}^{\ell}$ is proper and $\varphi(M)\subset B_N(r)\times\mathbb{R}^{\ell}$, then
\[
\sup_M|H|\geq\frac{m-\ell}{m}C_b(r).
\]
As a consequence, it follows from here that a complete hypersurface of given constant mean curvature lying inside a closed circular cylinder
in Euclidean space cannot be proper if the circular bases is of sufficiently small radius (see Corollary 1 in \cite{ABD}). In
particular, there exists no complete minimal hypersurface properly immersed in $\mathbb{R}^n$ and having 2 (or more) bounded
coordinates, showing that any possible counterexample to a conjecture of Calabi on complete minimal hypersurfaces \cite{Ca} (see
also \cite{chern}) cannot be proper.

Our main result here deals with the sectional curvature of such submanifolds and it can be stated as follows.
\begin{theorem}
\label{thmMain}
Let $M^{m}$ and $N^{n-\ell}$  be complete Riemannian manifolds of dimension $m$ and $n-\ell$ respectively, with
$n+\ell\leq 2m-1$. Let $\varphi:M^m\rightarrow N^{n-\ell}\times\R^{\ell}$ be a $($cylindrically
bounded$)$ isometric immersion with  $\varphi(M)\subset B_N(r)\times\R^{\ell}$.
Assume that the radial sectional curvature $K_N^{\mathrm{rad}}$ along the radial geodesics issuing from $p$
satisfies $K_N^{\mathrm{rad}}\leq b$ in $B_N(r)$ and $0<r<\min\{\mathrm{inj}_N(p),\pi/2\sqrt{b}\}$, where we replace
$\pi/2\sqrt{b}$ by $+\infty$ if $b\leq 0$. If either
\begin{itemize}
\item[(i)] the scalar curvature of $M$ satisfies (\ref{eqScalar}), or
\item[(ii)] $\varphi$ is proper and
\begin{equation}
\sup_{\varphi^{-1}(B_{N}(r)\times \partial B_{\mathbb{R}^{\ell}}(t))}\Vert\alpha\Vert\leq \sigma(t),
\label{growth}
\end{equation}
\end{itemize}
where $\alpha$ is the second fundamental form of the immersion and $\sigma:[0,+\infty)\rightarrow\mathbb{R}$ is a
positive  function satisfying $\int_0^{+\infty}1/\sigma=+\infty$, then
\begin{equation}
\sup_{M}K_{M}\geq C_{b}^{2}(r)+\inf_{B_{N}(r)}K_{N}.
\label{eq-ABM}
\end{equation}
\end{theorem}

It is worth pointing out that the codimension restriction $n+\ell\leq 2m-1$ cannot be relaxed. Actually, together with the
bound  $m\leq n-1$, it implies that $n-\ell\geq 3$ and $m\geq \ell +2$. In particular, for $n=3$  we have that
$\ell =0$, and therefore $\varphi (M)\subset B_{N}(r)$. In fact, the flat cylinder
$\mathbb{S}^{1}(r')\times \mathbb{R}\subset  B_{\mathbb{R}^{2}}(r)\times \mathbb{R}$, for every
$0<r'<r$, shows that the restriction $2m-1\geq n+\ell$ is necessary.

On the other hand, our estimate \rf{eq-ABM} is sharp. Indeed, for  every $n\geq 3$  and $\ell\leq n-3$ we can consider
$M=\partial B_{\mathbb{N}^{n-\ell}(b)}(r')\times \mathbb{R}^{\ell}$ and take
$\varphi\colon M\hookrightarrow B_{\mathbb{N}^{n-\ell}(b)}(r)\times \mathbb{R}^{\ell}$ to be
the canonical isometric immersion, for every $0<r'<r$. Therefore $\sup_MK_M$ is  the constant sectional
curvature of the geodesic sphere $\partial B_{\mathbb{N}^{n-\ell}(b)}(r')\subset\mathbb{N}^{n-\ell}(b)$, which is given by
\[
\sup_M K_M=K_{\partial B_{\mathbb{N}^{n-\ell}(b)}(r')}=
\begin{cases}
{b}/{\sin^2(\sqrt{b}\,r')} & \mbox{\rm if $b>0$ and $0<r'<\pi/2\sqrt{b}$}\\
1/r'^2 & \mbox{\rm if $b=0$ and $r'>0$}\\
{-b}/\sinh^2(\sqrt{-b}\,r') & \mbox{\rm if $b<0$ and $r'>0$}.
\end{cases}
\]
In particular, observe that
\[
\sup_M K_M=K_{\partial B_{\mathbb{N}^{n-\ell}(b)}(r')}=C_b^2(r')+b.
\]
Since in this case $K_{\mathbb{N}^{n-\ell}(b)}=b$, then for every $0<r'<r$ we have
\[
\sup_MK_M=C_b^2(r')+b\geq C_b^2(r)+\inf_{B_{\mathbb{N}^{n-\ell}(b)}(r)}K_{\mathbb{N}^{n-\ell}(b)},
\]
which shows that our estimate (\ref{eq-ABM}) is sharp.

Let $\varphi \colon M^{m}\to N^{n-\ell}\times \mathbb{R}^{\ell}$ be an isometric immersion of a compact Riemannian
$m$-manifold $M^{m}$, and let $\pi_{N}\colon N^{n-\ell}\times \mathbb{R}^{\ell}\to N^{n-\ell}$ be the
projection onto the first factor. Denote by $R_{M}$ the radius of the smallest ball of $N$ containing $\pi_{N}(\varphi (M))$. We will refer
to $R_M$ as the extrinsic radius of the immersion. As a consequence of Theorem\rl{thmMain}, we have the following versions of
the Extrinsic Radius Theorem of Jorge-Xavier \cite{jorge-xavier} (see also \cite[Theorem 1.3]{bessa-montenegro4}).
\begin{corollary}
\label{corollary-new}
Let $\varphi \colon M^{m}\to N^{n-\ell}\times \mathbb{R}^{\ell}$ be an isometric immersion of a compact Riemannian
$m$-manifold $M^{m}$ into the product $N^{n-\ell}\times\R^{\ell}$ with  $ n+\ell\leq 2m-1$, where $N^{n-\ell}$ is a complete
Riemannian manifold with a pole and radial sectional curvature $K_N^{\mathrm{rad}}\leq b\leq 0$.
Then, the extrinsic radius satisfies
\begin{equation}
R_{M}\geq C_{b}^{-1}\left(\sqrt{\sup K_{M}-\inf K_{N}}\right)
\end{equation}
In particular, if $N=\mathbb{R}^{n-\ell}$ we have that
\begin{equation}
R_{M}\geq \frac{1}{\sqrt{\sup K_{M}}}\cdot
\end{equation}
\end{corollary}

\begin{corollary}
Let $\varphi\colon M^{m}\to \mathbb{S}^{n-\ell}\times \mathbb{R}^{\ell}$ be an isometric immersion of a compact Riemannian
$m$-manifold $M^{m}$ with  $n+\ell\leq 2m-1$. If $\sup K_{M}\leq 1$ then
\begin{equation}
R_{M}\geq \pi/2.
\end{equation}
\end{corollary}

On the other hand, it is important to remark that for  hypersurfaces, the growth rate
of the norm of the second fundamental form in (\ref{growth}) can be improved as follows.
\begin{theorem}
\label{thMain2}
Let $\varphi:M^{n-1}\rightarrow N^{n-\ell}\times\R^{\ell}$ be a   properly immersed hypersurface
with $\varphi(M)\subset B_N(r)\times\R^{\ell}$, $n-\ell\geq 3$. Suppose that $N$ satisfies the assumptions on the radial sectional curvatures
as in Theorem\rl{thmMain} and the second fundamental form $\alpha$ satisfies
\begin{equation}
\label{growth2}
\sup_{\varphi^{-1}(B_{N}(r)\times \partial B_{\mathbb{R}^{\ell}}(t))} \Vert \alpha \Vert \leq  \sigma^{2}(t),
\end{equation}where $\sigma\colon [0, \infty)\to \mathbb{R}$ is a
positive  function satisfying
\[
\int_0^{+\infty}\frac{1}{\sigma}=+\infty, \quad {and} \quad
\limsup_{t\rightarrow+\infty}\displaystyle\frac{1}{\sigma(t)}<+\infty.
\]
Then
\[
\sup_{M}K_{M}\geq C_{b}^{2}(r)+\inf_{B_{N}(r)}K_{N}.
\]
\end{theorem}

\begin{remark}
It should be remarked  that Hasanis and Koutroufiotis \cite{hasanis-koutroufiotis} established similar sectional curvature estimates for cylindrically bounded submanifolds,  with scalar curvature bounded below, of the Euclidean space $\mathbb{R}^{n}$.
In a slightly  more general situation, F. Gim\'enez \cite{gimenez} established sectional curvature estimates for submanifolds
with scalar curvature bounded below immersed in a tubular neighborhood of  certain, ($P$-submanifolds),  embedded submanifolds
of Hadamard manifolds.   Our main results  besides extending  Hasanis-Koutrofiotis results to a larger class of submanifolds
can be easily adapted to reproduce Gim\'enez's result settings.
\end{remark}

\section{Preliminaries}

Our main tool  to build the  proof of Theorem \ref{thmMain} is the following (and important) version  of the
Omori-Yau Maximum Principle for the Hessian due to Pigola, Rigoli and Setti \cite[Theorem 1.9]{pigola-rigoli-setti}. Chronologically, the Omori-Yau Maximum Principle can be traced  in a series of papers, starting with
Omori \cite{omori} proving the maximum principle at infinity for the class of complete Riemannian manifolds with sectional curvature bounded from below, followed by new insights given by Cheng and Yau \cite{cheng-yau},
\cite{yau} extending to the class of complete Riemannian manifolds with Ricci curvature bounded from below and  extensions due to C. Dias \cite{dias} (extended to the class of sectional curvature with quadratic decay)  and
Chen-Xin \cite{chen-xin} (extended to the class of Ricci curvature with quadratic decay).
Finally, as observed by Pigola, Rigoli and Setti in \cite{pigola-rigoli-setti}, the validity of the Omori-Yau Maximum Principle does not depend
on curvature bounds as much as one would expect. Actually, a condition to guarantee the validity of it can be
expressed in a function theoretic form by the newest elaboration stated below.
\begin{theorem}\cite[Theorem 1.9]{pigola-rigoli-setti}
\label{TheoremPRS}
Let $M^m$ be a Riemannian manifold and assume that there exists a non-negative $C^2$-function $\psi$ satisfying the following requirements:
\begin{itemize}
\item[(a.1)] $\psi$ is proper, that is, $\psi(x)\rightarrow+\infty$ as $x\rightarrow\infty$;
\item[(a.2)] there exists a positive constant $A>0$ such that $|\grad\psi|\leq A\sqrt{\psi}$ outside a compact subset of $M$;
\item[(a.3)] there exists a positive constant $B>0$ such that
\[
\Hess\psi\leq B\sqrt{\psi G\left(\sqrt{\psi}\right)}\langle,\rangle
\]
(in the sense of quadratic forms) outside a compact subset of $M$,
\end{itemize}
where $G$ is a smooth function on $[0,+\infty)$ satisfying:
\begin{equation}\begin{array}{llll}
(i)\,\, G(0)>0, &  & & (ii) \,\, G'(t)\geq 0 \,\,{\it on}\,\,[0, +\infty),\vspace{1ex}\\
(iii)\,\,1/ \sqrt{G(t)}\not\in L^{1}(0,+\infty), &  & & (iv) \,\, \limsup_{t\to +\infty}\displaystyle{\frac{tG(\sqrt{t})}{G(t)}}<+\infty.
\end{array}\label{eqG}
\end{equation}
Then, given a function $u\in C^2(M)$ with  $u^{\ast}=\sup_{M}u <+\infty$ there exists a sequence
$\{x_{k}\}_{k\in\mathbb{N}}\subset M$ such that
\beq
\label{OY}
\textrm{$($i$)$} \,\,\, u(x_k)>u^*-\frac{1}{k}, \,\,\, \textrm{$($ii$)$} \,\,\,
|\grad u(x_k)|<\frac{1}{k},
\textrm{ and } \textrm{$($iii$)$} \,\,\, \Hess u(x_k)<\frac{1}{k}\langle,\rangle
\eeq
in the sense of quadratic forms, that is,
\[
\Hess u(x_k)(v,v)<\frac{1}{k}|v|^2 \,\, \text{ for all $v\in T_{x_k}M$}.
\]
If, instead of (a.3), one replaces it by
\begin{itemize}
\item[(a.3)'] there exists a positive constant $B>0$ such that
$\Delta\gamma\leq B\sqrt{\gamma G\left(\sqrt{\gamma}\right)}$ outside a compact subset of $M$,
\end{itemize}
then one can weaken conclusion (iii) in \rf{OY} to
\[
\textrm{$($iii$)'$} \,\,\, \Delta u(x_k)<\frac{1}{k}.
\]
\end{theorem}
\begin{remark}
It is worth pointing out that although in the statement (and in the proof) of Theorem \rl{TheoremPRS} the manifold $M$ is not
required to be geodesically complete, the two assumptions (a.1) and (a.2) imply it. For the details, see \cite[page 10]{pigola-rigoli-setti}.
\end{remark}

For that reason, and following the terminology introduced by Pigola, Rigoli and Setti in \cite{pigola-rigoli-setti}, the Omori-Yau maximum principle for the Hessian is said to hold on a (not necessarily complete)
Riemannian manifold $M$ if, for any smooth funtion $u\in C^2(M)$ with $u^*=\sup_M u<+\infty$ there exists a sequence of points $\{x_k\}_{k\in\mathbb{N}}\subset M$ satisfying
\[
\textrm{$($i$)$} \,\,\, u(x_k)>u^*-\frac{1}{k}, \,\,\, \textrm{$($ii$)$} \,\,\,
|\grad u(x_k)|<\frac{1}{k},
\textrm{ and } \textrm{$($iii$)$} \,\,\, \Hess u(x_k)<\frac{1}{k}\langle,\rangle
\]
for each $k\in\mathbb{N}$. Equivalently, for any smooth funtion $u\in C^2(M)$ with $u_*=\inf_M u>-\infty$ there exists a sequence of points
$\{x_k\}_{k\in\mathbb{N}}\subset M$ with the properties
\[
\textrm{$($i$)$} \,\,\, u(x_k)<u_*+\frac{1}{k}, \,\,\, \textrm{$($ii$)$} \,\,\,
|\grad u(x_k)|<\frac{1}{k},
\textrm{ and } \textrm{$($iii$)$} \,\,\, \Hess u(x_k)>-\frac{1}{k}\langle,\rangle
\]
for each $k\in\mathbb{N}$. In the same way, the Omori-Yau maximum principle is said to hold on a (not necessarily complete)
Riemannian manifold $M$ if, for any smooth funtion $u\in C^2(M)$ with $u^*=\sup_M u<+\infty$ there exists a sequence of points $\{x_k\}_{k\in\mathbb{N}}\subset M$ satisfying
\[
\textrm{$($i$)$} \,\,\, u(x_k)>u^*-\frac{1}{k}, \,\,\, \textrm{$($ii$)$} \,\,\,
|\grad u(x_k)|<\frac{1}{k},
\textrm{ and } \textrm{$($iii$)'$} \,\,\, \Delta u(x_k)<\frac{1}{k}
\]
for each $k\in\mathbb{N}$.

The function theoretic approach to the Omori-Yau maximum principle given in Theorem\rl{TheoremPRS} allows one to apply it in
different situations, where the choice of $\psi$ and $G$ are suggested by the geometric setting. For instance, one has the following
consequence (see \cite[Example 1.13]{pigola-rigoli-setti}).
\begin{corollary}
\label{example1.13}
Let $M$ be a complete, non-compact, Riemannian manifold, let $o\in M$ be a reference point and denote by $\varrho_M(x)$ de
Riemannian distance
function from $o$. Assume that the radial sectional curvature of $M$, that is, the sectional curvature of the 2-planes containing
$\grad \varrho_M$, satisfies
\[
K^{\text{rad}}_M\geq -G(r),
\]
where $G$ is a smooth function on $[0,+\infty)$ which we assume to be even at the origin, that is, $G^{(2k+1)}(0)=0$ for
$k=0,1,2,\ldots$,
and satisfying the conditions listed in Theorem\rl{TheoremPRS}. Then the Omori-Yau maximum principle for the Hessian holds on $M$.
\end{corollary}

We will need two  more results. The first is known as Otsuki's lemma (for a proof see, for instance,
\cite[Page 28]{kobayashi-nomizu} or \cite[Lemma 1.16]{pigola-rigoli-setti}).
\begin{lemma}\label{lemmaOtsuki}
Let $\beta\colon\mathbb{R}^{k}\times \mathbb{R}^{k}\to \mathbb{R}^{q}$, $q\leq k-1$,
be a bilinear symmetric form satisfying $\beta(X,X)\neq 0$ for $X\neq 0$.
Then there exist linearly independent vectors $X$, $Y$ such that
\beq
\beta(X,X)=\beta(Y,Y) \quad \text{and} \quad  \beta(X,Y)=0.
\eeq
\end{lemma}

And the second is the well known Hessian Comparison Theorem, see \cite{schoen-yau}.
\begin{theorem}[Hessian Comparison Theorem]
\label{thm2}
Let $M$ be a  Riemannian manifold  and
$x_0,x_1 \in M$ be such that there is a minimizing unit speed geodesic $\gamma$  joining $x_{0}$ and $x_{1}$ and let $\varrho(x)=\mathrm{dist}(x_0,x)$
be the  distance function  to $x_{0}$. Let $K_{\gamma}\leq b$ be the radial sectional curvatures of $M$ along $\gamma$.
If $b>0$ assume $\varrho(x_{1})<\pi/2\sqrt{b}$. Then, we have
\[
\Hess\varrho(x) (\gamma',\gamma')=0
\]
and
\[
\Hess\varrho(x)(X,X)\geq C_{b}(\varrho(x))|X|^2
\]
where $X\in T_{x}M$ is perpendicular to $\gamma'(\varrho(x))$.
\end{theorem}

\section{Proof of Theorem \ref{thmMain}}

Let $\varphi\colon M^m\to\overline{M}^n$ be an isometric immersion between Riemannian manifolds. Given a function
$g\in C^\infty(\overline{M})$ we set $f=g\circ\varphi\in C^\infty(M)$.
Since
\[
\langle\grad^Mf,X\rangle=\langle\grad^{\overline{M}}g, X\rangle
\]
for every vector field $X\in TM$, we obtain
\[
\grad^{\overline{M}}g=\grad^Mf +(\grad^{\overline{M}}g)^{\perp}
\]
according to the decomposition $T\overline{M}=TM\oplus T^\perp M$. An easy computation using the Gauss formula gives the
well-known relation (see e.g. \cite{jorge-koutrofiotis})
\begin{equation}
\label{eqBF2}
\Hess^Mf(X,Y)= \Hess^{\overline{M}}g(X,Y) +\langle\grad^{\overline{M}}g,\alpha(X,Y)\rangle
\end{equation}
for all vector fields $X,Y\in TM$, where $\alpha$ stands for the second fundamental form of
$\varphi$. In particular, taking traces  with respect to an orthonormal frame
$\{ e_{1},\ldots, e_{m}\}$ in $TM$ yields
\[
\Delta_Mf =\sum_{i=1}^{m}\Hess^{\overline{M}}g(e_i,e_i)+ \langle\grad^{\overline{M}}g,\stackrel{\to}{H}\rangle.
\]
where $\stackrel{\to}{H}=\sum_{i=1}^{m}\alpha(e_i,e_i)$.

\subsection{Proof of Theorem \ref{thmMain}, item (i)}
Let $g:N\times\mathbb{R}^\ell\rightarrow\mathbb{R}$ be given $g(z,y)=\phi_b(\varrho_N(z))$, where
\[
\phi_b(t)=
\begin{cases}
1-\cos(\sqrt{b}\,t) &\text{if $b>0$ and $0<t<\pi/2\sqrt{b}$}\\
t^2 &\text{if $b=0$ and $t>0$}\\
\cosh(\sqrt{-b}\,t) &\text{if $b<0$ and $t>0$},
\end{cases}
\]
and $\varrho_N(z)=\mathrm{dist}_N(p,z)$. Consider $f:M\rightarrow\R$ the function $f=g\circ\varphi$,
and let $\pi_N:N\times\R^{\ell}\rightarrow N$ be the projection on the factor $N$.
Since $\pi_N(\varphi(M))\subset B_N(r)$, we have that $f^*=\sup_Mf\leq\phi_b(r)<+\infty$. The idea of the proof is similar to the idea of Jorge-Koutrofiotis in \cite{jorge-koutrofiotis}. We will
need to apply the Omori-Yau maximum principle for the Hessian to the function $f$ in order to control the second fundamental form of the immersion restricted to a certain subspaces and apply Otsuki's Lemma in the  estimate
of the sectional curvature.

To see that the Omori-Yau maximum principle for the Hessian holds on $M$ we may suppose that $\sup K_{M}<\infty$.
Otherwise the estimate in (\ref{eq-ABM}) is trivially satisfied. In that case, since the scalar curvature is an average of
sectional curvatures it follows from (\ref{eqScalar}) that
\begin{equation}
K^{\mathrm{rad}}_{M}(x)\geq -\widehat{B}^2\varrho^{2}_{M}(x)\cdot\prod_{j=1}^{k}\Big(\log^{(j)}(\varrho_{M}(x))\Big)^2,\,\, \varrho_{M}(x)\gg 1,
\label{eqKrad}
\end{equation}
for a positive constant $\widehat{B}>0$, where $K^{\mathrm{rad}}_{M}$  denotes the radial sectional curvature of $M$.
According to Corollary \ref{example1.13}, this curvature decay suffices to conclude that the Omori-Yau Maximum Principle for
the Hessian holds on $M$.
Therefore, there exists a sequence of points $\{x_{k}\}$ in $M$ such that
\begin{equation}
\label{eqluis21}
\textrm{(i)} \,\,\, f(x_k)>f^*-\frac{1}{k}, \,\,\, \textrm{(ii)} \,\,\,
|\grad f(x_k)|<\frac{1}{k},
\textrm{ and } \textrm{(iii)} \,\,\, \Hess f(x_k)<\frac{1}{k}\langle,\rangle.
\end{equation}

Since $f(x)=g(\varphi(x))=\phi_b(\varrho_N(z))$, where $z=z(x)=\pi_N(\varphi(x))$, then
\beq
\label{eqPacelli}
\gradien^{N\times\mathbb{R}^\ell}g(\varphi(x))=\grad f(x)+(\gradien^{N\times\mathbb{R}^\ell}g(\varphi(x)))^\perp,
\eeq
where
\[
\gradien^{N\times\mathbb{R}^\ell}g(z,y)=\phi'_b(\varrho_N(z))\gradien^{N}\varrho_N(z).
\]
It then follows from \rf{eqBF2} that
\begin{eqnarray}
\label{eqluis1}
\nonumber \Hess f(x)(X,X) & = & \hessian^N(\phi_b\circ\varrho_N)(z)(\pi_{TN}X,\pi_{TN}X)\\
{} & {} & {} \\
\nonumber {} & {} &  +\langle\gradien^N(\phi_b\circ\varrho_N)(z),\alpha_x(X,X)\rangle
\end{eqnarray}
for all vector field $X\in T_xM$, where $\pi_{TN}$ denotes the orthogonal projection of $TM$ onto $TN$.
Observe also that
\begin{eqnarray}
\label{eqluis2} \hessian^N(\phi_b\circ\varrho_N)(z)(\pi_{TN}X,\pi_{TN}X) & = &
 \phi''_b(s)\left(\frac{\partial\varrho_N}{\partial X}\right)^2\nonumber \\
 {} & {} &\\
 {} & + & \phi'_b(s)\hessian^N\varrho_N(z)(\pi_{TN}X,\pi_{TN}X)\nonumber
\end{eqnarray}
where $s=s(x)=\varrho_N(z)$ and
\[
\frac{\partial\varrho_N}{\partial X}=\langle\gradien^N\varrho_N,X\rangle.
\]
Taking into account that the function $\phi_b(t)$ satisfies the differential equation
\[
\phi''_b(t)-\phi'_b(t)C_b(t)=0,
\]
it follows from (\ref{eqluis1}) and (\ref{eqluis2}) that
\begin{eqnarray}
\label{eqluis3}
\nonumber \Hess f(x)(X,X) & = &
\phi'_b(s)\left(C_b(s)\left(\frac{\partial\varrho_N}{\partial X}\right)^2+\langle\gradien^N\varrho_N(z),\alpha_x(X,X)\rangle\right) \\
{} &  & \\
{} & + & \phi'_b(s)\hessian^N\varrho_N(z)(\pi_{TN}X,\pi_{TN}X)\nonumber
\end{eqnarray}

Since $m\geq \ell+2$, we have for each $x\in M$ a subspace $V_x\subset T_xM\subset T_{\varphi(x)}(N\times\R^{\ell})$
such that $V_x\perp T\R^\ell$ and $\mathrm{dim} V_x\geq m-\ell\geq 2$. Choose
$\{\partial/\partial\varrho_N,\partial/\partial\theta_{2},\ldots,\partial/\partial\theta_{n-\ell}\}$ orthonormal
polar coordinates for $TN$. Then, for every $X\in V_x$ we have $\pi_{T\mathbb{R}^\ell}X=0$ and
\[
X=a_1^X\frac{\partial}{\partial\varrho_N}+\sum_{j=2}^{n-\ell}a_j^X\frac{\partial}{\partial\theta_j},
\]
where $a_1^X=\partial\varrho_N/\partial X$. Therefore, using Theorem \ref{thm2}, we have that for every $X\in V_x$
\begin{eqnarray*}
\hessian^N\varrho_N(z)(\pi_{TN}X,\pi_{TN}X) & = &
\sum_{j=2}^{n-\ell}(a_j^X)^2\hessian^N\varrho_N(z)(\frac{\partial}{\partial\theta_j},\frac{\partial}{\partial\theta_j})\\
{} & \geq & C_b(s)\left(\langle , \rangle-d\varrho_N\otimes d\varrho_N\right)(X,X)\\
{} & = & C_b(s)(|X|^2-(a_1^X)^2),
\end{eqnarray*}
since $\pi_{TN}X=X$, so that by (\ref{eqluis3}) we have
\begin{eqnarray}
\label{eqluis22}
\Hess f(x)(X,X) & \geq &
\phi'_b(s)\left(C_b(s)|X|^2+\langle\gradien^N\varrho_N(z),\alpha_{x}(X,X)\rangle\right)\nonumber \\
&&\\
{} & \geq &
\phi'_b(s)\left(C_b(s)|X|^2-|\alpha_{x}(X,X)|\right)\nonumber
\end{eqnarray}
for every $X\in V_x$. From here and (\ref{eqluis21}) we obtain that
\[
\frac{1}{k}|X|^2\geq \Hess f(x_k)(X,X)\geq
\phi'_b(s_k)\left(C_b(s_k)|X|^2-|\alpha_{x_k}(X,X)|\right)
\]
for every $x_k$ and every $X\in V_{x_k}$, where $z_k=\pi_N(\varphi(x_k)$ and $s_k=s(x_k)=\varrho_N(z_k)$. Hence
\[
|\alpha_{x_k}(X,X)|\geq \left(C_b(s_k)-\frac{1}{k\phi'_b(s_k)}\right)|X|^2
\]
with
\[
C_b(s_k)-\frac{1}{k\phi'_b(s_k)}>0
\]
for $k$ sufficiently large.

Now consider $\beta_{x_k}\colon V_{x_{k}}\times V_{x_{k}}\to \mathbb{R}^{n-m}$ the restriction of the second fundamental
form $\alpha_{x_k}$  to $V_{x_{k}}$. We have that
\[
n-m\leq m-\ell-1\leq\mathrm{dim}V_{x_k}-1
\]
since  $2m-1\geq n+\ell$, and
therefore we may apply Lemma \ref{lemmaOtsuki} to $\beta_{x_k}$. The conclusion is that there are linearly independent vectors
$X_k,Y_k\in V_{x_{k}}$ such that
\[
\alpha(X_k,X_k)=\alpha (Y_k,Y_k) \quad \text{and} \quad \alpha (X_k,Y_k)=0
\]
and $|X_k|\geq |Y_k|\geq 1$. We will now compare the sectional curvature $K_{M}(X_k,Y_k)$ in $M$ of the plane spanned by
$X_k$ and $Y_k$ with the sectional curvature $K_{N\times\mathbb{R}^\ell}(X_k,Y_k)$ in $N\times\mathbb{R}^\ell$ of the same plane.
Since $X_k,Y_k\in V_{x_k}\perp T\mathbb{R}^\ell$, then
\[
K_{N\times\mathbb{R}^\ell}(X_k,Y_k)=K_{N}(X_k,Y_k).
\]
Then, by the Gauss equation we have that
\begin{eqnarray*}
K_{M}(X_k,Y_k)-K_{N}(X_k,Y_k) & = & K_{M}(X_k,Y_k)-K_{N\times\mathbb{R}^\ell}(X_k,Y_k)\\
{} & = & \frac{\langle \alpha (X_k,X_k), \alpha (Y_k,Y_k)\rangle -\vert \alpha (X_k,Y_k)\vert^{2}}{\vert X_k\vert^{2}\vert Y_k\vert^{2}-\langle X_k, Y_k\rangle^{2} } \\
& \geq & \frac{\vert  \alpha (X_k,X_k)\vert^2}{\vert X_k\vert^{2}\vert Y_k\vert^{2} }\geq
\left(\frac{\vert \alpha (X_k,X_k)\vert}{|X_k|^2}\right)^2  \\
&\geq & \left(C_{b}(s_{k})-\frac{1}{k\phi'_b(s_k)}\right)^{2}\!\!\!\cdot
\end{eqnarray*}
Thus
\begin{equation}
\sup_{M}K_{M}-\inf_{B_{N}(r)}K_{N}\geq \left(C_{b}(s_{k}) -\frac{1}{k\,\phi_{b}'(s_{k})}\right)^{2}\!\!\!\cdot
\end{equation}
Observe that $f^*=\phi_b(s^*)$, where $s^*=\sup_M s$ and $s_k\to s^*\leq r$. Therefore, letting $k\to \infty $ we have
that
\[
\sup_{M}K_{M}-\inf_{B_{N}(r)}K_{N}\geq C_{b}^{2}(s^*)\geq C_b^2(r).
\]
This finishes the proof of item (i) in Theorem \ref{thmMain}.

\subsection{Proof of Theorem \ref{thmMain}, item (ii)}
In this case, we cannot apply directly Theorem\rl{TheoremPRS}, but we may apply parts of its proof. Consider again $f:M\rightarrow\mathbb{R}$
the function given by $f(x)=\phi_b(\varrho_N(z(x)))$, with $z(x)=\pi_N(\varphi(x))$. Since $\pi_N(\varphi(M))\subset B_N(r)$, we
have that $f^*=\sup_Mf\leq\phi_b(r)$.
Let $\psi:M\to [0,+\infty)$ be given by
\[
\psi(x)=\exp\left(\int_0^{|y(x)|}\frac{ds}{\sigma(s)}\right),
\]
where $y(x)=\pi_{\mathbb{R}^\ell}(\varphi(x))$. Since $\varphi$ is proper and $\pi_N(\varphi(M))\subset B_N(r)$, then the
function $|y(x)|$ satisfies $|y(x)|\rightarrow+\infty$ as $x\rightarrow\infty$. By hypothesis we have that
$\int_0^{+\infty}1/\sigma(s)ds=+\infty$, so that  $\psi(x)\to+\infty$ as $x\to\infty$.

Following the ideas of Pigola, Rigoli and Setti  in the proof of \cite[Theorem 1.9]{pigola-rigoli-setti}, we let
$x_0\in M$ with $\pi_{N}(\varphi(x_0))\neq p$ and set
$$
f_{k}(x)=\frac{f(x)-f(x_0)+1}{\psi(x)^{1/k}}.
$$
Thus $f_{k}(x_0)>0$ and since $f^*\leq\phi_{b}(r)<+\infty$ and $\psi(x)\to+\infty $ as $x\to \infty$, we have that
$\limsup_{x\to\infty}f_{k}(x)\leq 0$. Hence $f_{k}$ attains a positive absolute maximum at a point $x_{k}\in M.$ This
procedure yields a sequence $\{x_{k}\}$ such that (passing to a subsequence if necessary) $f(x_{k})$ converges to
$f^*$ (see page 8 of \cite{pigola-rigoli-setti}).

Suppose first that $x_{k}\to\infty$ as $k\to+\infty$. Since $f_{k}$ attains a positive maximum at $x_{k}$ we have
$\grad f_{k}(x_{k})=0$ and $\Hess f_{k}(x_{k})(X,X)\leq 0$ for every $X\in T_{x_k}M$. This yields
\begin{equation}
\label{eqB}
\grad f(x_{k})= \frac{f(x_{k})-f(x_0)+1}{k\psi (x_{k})}\grad \psi (x_{k})
\end{equation}
and
\begin{eqnarray}
\label{eqA}
\nonumber \Hess f (x_{k}) & \leq & \frac{f(x_{k})-f(x_0)+1}{k\psi(x_k)}
\left(\Hess\psi(x_{k})+\left(\frac{1}{k}-1\right)\frac{1}{\psi(x_{k})}d\psi\otimes d\psi \right)\\
&&\\
{} & \leq & \frac{f(x_{k})-f(x_0)+1}{k\psi(x_k)}\Hess\psi(x_{k}).\nonumber
\end{eqnarray}

Since $\psi(x)=\zeta(y)$ where $y=y(x)$ and $\zeta(y)=\exp(\int_0^{|y|}ds/\sigma(s))$, $y\in\mathbb{R}^{\ell}$, from
\rf{eqBF2} we have that
\begin{equation}
\label{eqluis23}
\Hess\psi(x)(X,X)=\hessian^{\mathbb{R}^\ell}
\zeta(y)(\pi_{T\mathbb{R}^\ell}X,\pi_{T\mathbb{R}^\ell}X)
+\langle\gradien^{\mathbb{R}^\ell}\zeta(y),\alpha_x(X,X)\rangle
\end{equation}
for all vectors $X\in T_xM$, where $\pi_{T\mathbb{R}^\ell}$ denotes the orthogonal projection of $TM$ onto
$T\mathbb{R}^\ell$. Observe also that
\begin{equation*}
\gradien^{\mathbb{R}^\ell}\zeta(y)=\frac{\zeta(y)}{\sigma(|y|)}\gradien^{\mathbb{R}^\ell}|y|,
\end{equation*}
and then
\begin{equation}
\label{eqC}
\grad\psi(x)=\frac{\psi(x)}{\sigma(|y|)}\gradien^{\mathbb{R}^\ell}|y|.
\end{equation}
Thus, for every $X\in T_xM$ such that $\pi_{T\mathbb{R}^\ell}X=0$ it follows from (\ref{eqluis23}) that
\begin{equation*}
\Hess\psi(x)(X,X)=\frac{\psi(x)}{\sigma(|y(x)|)}\langle\gradien^{\mathbb{R}^\ell}|y|,\alpha_x(X,X)\rangle\leq
\frac{\psi(x)}{\sigma(|y(x)|)}|\alpha_x(X,X)|.
\end{equation*}
Therefore, by (\ref{growth}) we obtain that
\begin{equation}
\label{eqluis24}
\frac{1}{\psi(x)}\Hess\psi(x)(X,X)\leq\frac{|\alpha_x(X,X)|}{\sigma(|y(x)|)}\leq |X|^2
\end{equation}
for every $X\in T_xM$ with $\pi_{T\mathbb{R}^\ell}X=0$.

As in the proof of item (i), since $m\geq \ell+2$, we may choose for each ${x_k}\in M$ a subspace
$V_{x_k}\subset T_{x_k}M$ with $\mathrm{dim} V_{x_k}\geq m-\ell\geq 2$ and such that
$V_{x_k}\perp T\R^\ell$. Then, $\pi_{T\mathbb{R}^\ell}X=0$ for every $X\in V_{x_k}$, and from
(\ref{eqA}) and (\ref{eqluis24}) we get that
\begin{equation}
\Hess f (x_{k})(X,X)\leq\frac{f(x_{k})-f(x_0)+1}{k\psi(x_k)}\Hess\psi(x_{k})(X,X)\leq\frac{\phi_b(r)+1}{k}|X|^2,
\end{equation}
for every $X\in V_x$. Moreover, using Theorem\rl{thm2}, we also have here that
\begin{equation}
\label{eqluis22bis}
\Hess f(x)(X,X)\geq \phi'_b(s)\left(C_b(s)|X|^2-|\alpha_{x}(X,X)|\right)
\end{equation}
for every $X\in V_x$, since $\pi_{TN}X=X$. Therefore, we obtain that
\[
\frac{\phi_b(r)+1}{k}|X|^2\geq \Hess f(x_k)(X,X)\geq
\phi'_b(s_k)\left(C_b(s_k)|X|^2-|\alpha_{x_k}(X,X)|\right)
\]
for every $x_k$ and every $X\in V_{x_k}$, where $z_k=\pi_N(\varphi(x_k))$ and $s_k=s(x_k)=\varrho_N(z_k)$. Hence
\[
|\alpha_{x_k}(X,X)|\geq \left(C_b(s_k)-\frac{\phi_b(r)+1}{k\phi'_b(s_k)}\right)|X|^2
\]
with
\[
C_b(s_k)-\frac{\phi_b(r)+1}{k\phi'_b(s_k)}>0
\]
for $k$ sufficiently large. Reasoning now as in the last part of the proof of item (i), there exist linearly independent vectors
$X_k,Y_k\in V_{x_k}$ such that, by Gauss equation,
\[
K_{M}(X_k,Y_k)-K_{N}(X_k,Y_k)=\left(\frac{\vert \alpha (X_k,X_k)\vert}{|X_k|^2}\right)^2
\geq \left(C_{b}(s_{k})-\frac{\phi_b(r)+1}{k\phi'_b(s_k)}\right)^{2}.
\]
We obtain from here that
\begin{equation}
\sup_{M}K_{M}-\inf_{B_{N}(r)}K_{N}\geq \left(C_{b}(s_{k}) -\frac{\phi_b(r)+1}{k\,\phi_{b}'(s_{k})}\right)^{2},
\end{equation}
and letting $k\to \infty $ we conclude that
$$\sup_{M}K_{M}-\inf_{B_{N}(r)}K_{N}\geq C_{b}^{2}(s^*)\geq C_b^2(r),$$
where $s^*=\sup_M s$,  $f^*=\phi_b(s^*)$ and $s_k\to s^*\leq r$.

To finish the proof of item (ii), we need to consider the case where the sequence $\{x_{k}\}\subset M$ remains in
a compact set. In that case, passing to a subsequence if necessary, we may assume that  $x_{k}\to x_{\infty}\in M$ and $f$
attains its absolute maximum at $x_{\infty}$. Thus $\Hess f (x_{\infty})(X,X)\leq 0$ for all $X\in T_{x_{\infty}}M$. In
particular, it follows from (\ref{eqluis22bis}) that for every $X\in V_{x_{\infty}}$
\[
0\geq \Hess f (x_{\infty})(X,X)\geq
\phi'_b(s_\infty)\left(C_b(s_\infty)|X|^2-|\alpha_{x_\infty}(X,X)|\right),
\]
where $s_\infty=\varrho_N(\pi_N(\varphi(x_\infty)))$.
Therefore
$$\vert \alpha_{x_\infty}(X,X)\vert\geq C_{b}(s_{\infty})|X|^2.$$
By applying Lemma \ref{lemmaOtsuki} to $\beta_{x_\infty}\colon V_{x_{\infty}}\times V_{x_{\infty}}\to \mathbb{R}^{n-m}$, the
restriction of the second fundamental form $\alpha_{x_\infty}$  to $V_{x_{\infty}}$, and reasoning now again as in the last part of
the proof of item (i), we have that there exist linearly independent vectors
$X_\infty,Y_\infty\in V_{x_\infty}$ such that, by Gauss equation,
\[
K_{M}(X_\infty,Y_\infty)-K_{N}(X_\infty,Y_\infty)=\left(\frac{\vert \alpha (X_\infty,X_\infty)\vert}{|X_\infty|^2}\right)^2
\geq C_{b}^2(s_{\infty}).
\]
Thus, we conclude from here that
\begin{equation}
\sup_{M}K_{M}-\inf_{B_{N}(r)}K_{N}\geq C_{b}^2(s_{\infty})\geq C_{b}^2(r).
\end{equation}
This finishes the proof of Theorem \ref{thmMain}.

\section{Proof of Theorem\rl{thMain2}} We proceed as in the proof of Theorem\rl{thmMain}, item (ii), to obtain a sequence
$\{x_{k}\}$ such that $f(x_{k})$ converges to $f^*$ and satisfying
\begin{equation}
\label{eqBbis}
\grad f(x_{k})= \frac{f(x_{k})-f(x_0)+1}{k\psi (x_{k})}\grad \psi (x_{k})
\end{equation}
and
\begin{equation}
\label{eqAbis}
\Hess f (x_{k})\leq\frac{f(x_{k})-f(x_0)+1}{k\psi(x_k)}\Hess\psi(x_{k}).
\end{equation}
Recall that (see \rf{eqC})
\begin{equation}
\label{eqCbis}
\grad\psi(x)=\frac{\psi(x)}{\sigma(|y|)}\gradien^{\mathbb{R}^\ell}|y|.
\end{equation}

Let us consider first the case where $x_{k}\to\infty$ as $k\to+\infty$. From \rf{eqBbis} and \rf{eqCbis} we know that
\[
|\grad f(x_k)|\leq\frac{(f^*+1)}{k}\frac{1}{\sigma(|y_k|)}\leq\frac{(\phi_b(r)+1)}{k}\frac{1}{\sigma(|y_k|)},
\]
where $y_k=y(x_k)$. Since $\varphi$ is proper and $\pi_N(\varphi(M))\subset B_N(r)$, then $|y_k|\rightarrow+\infty$
as $k\rightarrow+\infty$. Therefore, taking into account that $\limsup_{t\rightarrow+\infty}1/\sigma(t)<+\infty$ we obtain
from here that
\beq
\lim_{k\rightarrow+\infty}|\grad f(x_k)|=0.
\eeq
Observe that
\[
\gradien^{N\times\mathbb{R}^\ell}g(\varphi(x))=
\phi'_b(\varrho_N(z))\gradien^{N}\varrho_N(z)=
\grad f(x)+(\gradien^{N\times\mathbb{R}^\ell}g(\varphi(x)))^\perp,
\]
where $z=z(x)=\pi_N(\varphi(x))$. Therefore,
\beq
\label{eqPacellibis}
\phi'_b(s_k)^2=|\grad f(x_k)|^2+|(\gradien^{N\times\mathbb{R}^\ell}g(\varphi(x_k)))^\perp|^2,
\eeq
with $s_k=\varrho_N(z(x_k))$, and making $k\to\infty$ here we obtain that
\[
\lim_{k\rightarrow+\infty}|(\gradien^{N\times\mathbb{R}^\ell}g(\varphi(x_k)))^\perp|=\phi'_b(s^*)>0,
\]
which implies that $$(\gradien^{N\times\mathbb{R}^\ell}g(\varphi(x_k)))^\perp\neq 0$$ for $k$ sufficiently large.

As in the proof of Theorem\rl{thmMain}, since $m=n-1\geq \ell+2$, we may choose for each ${x_k}\in M$ a subspace
$V_{x_k}\subset T_{x_k}M$ with $\mathrm{dim} V_{x_k}\geq n-1-\ell\geq 2$ and such that
$V_{x_k}\perp T\R^\ell$.
Then, using Theorem\rl{thm2}, we also have here that
\beq
\label{eqA1}
\Hess f(x_k)(X,X)\geq \phi'_b(s_k)\left(C_b(s_k)|X|^2-|\alpha_{x_k}(X,X)|\right)
\eeq
for every $X\in V_{x_k}\leq T_{x_k}M$, since $\pi_{TN}X=X$.
On the other hand, we also know from \rf{eqAbis} that
\begin{eqnarray}
\label{eqA2}
\nonumber \Hess f(x_k)(X,X) & \leq & \frac{(\phi_b(r)+1)}{k}\frac{\Hess\psi(x_k)(X,X)}{\psi(x_k)}\\
{} & {} & {} \\
\nonumber {} & = & \frac{(\phi_b(r)+1)}{k}\frac{1}{\sigma(|y_k|)}\langle\gradien^{\mathbb{R}^\ell}|y|,\alpha_{x_k}(X,X)\rangle
\end{eqnarray}
for every $X\in T_{x_k}M$. Since $m=n-1$ and $(\gradien^{N\times\mathbb{R}^\ell}g(\varphi(x_k)))^\perp\neq 0$
(for $k$ large enough), then
\begin{equation}\label{eqPacelli2}
\alpha_{x_k}(X,X)=\lambda_{x_k}(X,X)(\gradien^{N\times\mathbb{R}^\ell}g(\varphi(x_k)))^\perp
\end{equation}
for a real function $\lambda$.

Observe now that
\begin{eqnarray*}
\langle\gradien^{\mathbb{R}^\ell}|y|,\alpha_{x_k}(X,X)\rangle & = & \lambda_{x_k}(X,X)
\langle\gradien^{\mathbb{R}^\ell}|y|,(\gradien^{N\times\mathbb{R}^\ell}g(\varphi(x_k)))^\perp\rangle \\
{} & = & \lambda_{x_k}(X,X)
\langle\gradien^{\mathbb{R}^\ell}|y|,\grad f(x_k)\rangle
\end{eqnarray*}
because of $\langle\gradien^{\mathbb{R}^\ell}|y|,\gradien^N\varrho_N\rangle=0$. Therefore,
\begin{eqnarray*}
\langle\gradien^{\mathbb{R}^\ell}|y|,\alpha_{x_k}(X,X)\rangle & \leq &  |\lambda_{x_k}(X,X)||\grad f(x_k)|\\
{} & \leq & |\lambda_{x_k}(X,X)|\frac{(\phi_b(r)+1)}{k}\frac{1}{\sigma(|y_k|)}.
\end{eqnarray*}
On the other hand, from our hypothesis \rf{growth2} we know that
\[
|\alpha_{x}(X,X)|\leq\sigma^2(|y(x)|)|X|^2,
\]
and from (\ref{eqPacellibis}) and (\ref{eqPacelli2}) we have that
\[
|\alpha_{x_k}(X,X)|=|\lambda_{x_k}(X,X)|\sqrt{\phi'_b(s_k)^2-|\grad f(x_k)|^2}\leq\sigma^2(|y_k|)|X|^2.
\]
That is,
\[
\frac{|\lambda_{x_k}(X,X)|}{\sigma(|y_k|)}\leq\frac{\sigma(|y_k|)|X|^2}{\sqrt{\phi'_b(s_k)^2-|\grad f(x_k)|^2}}.
\]
It follows from here that
\[
\langle\gradien^{\mathbb{R}^\ell}|y|,\alpha_{x_k}(X,X)\rangle \leq
\frac{(\phi_b(r)+1)}{k}\frac{\sigma(|y_k|)|X|^2}{\sqrt{\phi'_b(s_k)^2-|\grad f(x_k)|^2}}
\]
for every $X\in T_{x_k}M$, so that by \rf{eqA2} we get
\begin{equation}
\label{eqA3}
\Hess f(x_k)(X,X)\leq\frac{(\phi_b(r)+1)^2}{k^2}\frac{|X|^2}{\sqrt{\phi'_b(s_k)^2-|\grad f(x_k)|^2}}.
\end{equation}
Therefore, from \rf{eqA1} and \rf{eqA3} we have that
\beq
\label{eqA4}
\phi'_b(s_k)\left(C_b(s_k)|X|^2-|\alpha_{x_k}(X,X)|\right)\leq\frac{(\phi_b(r)+1)^2}{k^2}\frac{|X|^2}{\sqrt{\phi'_b(s_k)^2-|\grad f(x_k)|^2}}
\eeq
for every $X\in V_{x_k}$. Hence
\[
|\alpha_{x_k}(X,X)|\geq\left(C_b(s_k)-\frac{(\phi_b(r)+1)^2}{k^2\phi'_b(s_k)\sqrt{\phi'_b(s_k)^2-|\grad f(x_k)|^2}}\right)|X|^2,
\]
with
\[
\lim_{k\rightarrow+\infty}\left(C_b(s_k)-\frac{(\phi_b(r)+1)^2}{k^2\phi'_b(s_k)\sqrt{\phi'_b(s_k)^2-|\grad f(x_k)|^2}}\right)=
C_b(s^*)\geq C_b(r)>0,
\]
where $s^*=\sup_M s$,  $f^*=\phi_b(s^*)$ and $s_k\to s^*\leq r$.
Reasoning now as in the last part of the proof of item (i), there exist linearly independent vectors
$X_k,Y_k\in V_{x_k}$ such that, by Gauss equation,
\begin{eqnarray*}
K_{M}(X_k,Y_k)-K_{N}(X_k,Y_k) & = & \left(\frac{\vert \alpha (X_k,X_k)\vert}{|X_k|^2}\right)^2 \\
{} & \geq & \left(C_b(s_k)-\frac{(\phi_b(r)+1)^2}{k^2\phi'_b(s_k)\sqrt{\phi'_b(s_k)^2-|\grad f(x_k)|^2}}\right).
\end{eqnarray*}
We obtain from here that
\begin{equation}
\sup_{M}K_{M}-\inf_{B_{N}(r)}K_{N}\geq
\left(C_b(s_k)-\frac{(\phi_b(r)+1)^2}{k^2\phi'_b(s_k)\sqrt{\phi'_b(s_k)^2-|\grad f(x_k)|^2}}\right)^{2},
\end{equation}
and letting $k\to \infty $ we conclude that
$$\sup_{M}K_{M}-\inf_{B_{N}(r)}K_{N}\geq C_{b}^{2}(s^*)\geq C_b^2(r).$$

Finally, in the case where the sequence $\{x_k\}\subset M$ remains in a compact subset of $M$, and passing to a subsequence if
necessary, we may assume that  $x_{k}\to x_{\infty}\in M$ and $f$ attains its absolute maximum at $x_{\infty}$.
Thus, $\Hess f (x_{\infty})(X,X)\leq 0$ for all $X\in T_{x_{\infty}}M$.
Therefore, it follows again from Theorem\rl{thm2} that for every $X\in V_{x_{\infty}}$
\[
0\geq \Hess f (x_{\infty})(X,X)\geq
\phi'_b(s_\infty)\left(C_b(s_\infty)|X|^2-|\alpha_{x_\infty}(X,X)|\right),
\]
where $s_\infty=\varrho_N(\pi_N(\varphi(x_\infty)))$ and $V_{x_\infty}\subset T_{x_\infty}M$ is a subspace
with $\mathrm{dim} V_{x_\infty}\geq n-1-\ell\geq 2$ and such that $V_{x_\infty}\perp T\R^\ell$.
The proof now finishes as at the end of item (ii) in Theorem\rl{thmMain}.


\begin{thebibliography}{aa}



\bibitem{ABD} L.J. Al\'{\i}as, G.P. Bessa and M. Dajczer, \textit{The mean curvature of cylindrically bounded submanifolds}. Math. Ann.
{\bf 345} (2009), 367--376.

\bibitem{bessa-montenegro4} G. P. Bessa and J. F. Montenegro, \textit{On compact H-hypersurfaces of $N\times\R$}. Geom.
Dedicata. {\bf 127} (2007), 1--5.

\bibitem{Ca}E. Calabi, {\em Problems in Differential Geometry} (S. Kobayashi and J. Eells, Jr., eds.) Proc. of the United States-Japan Seminar in Differential Geometry, Kyoto, Japan, 1965, Nippon Hyoronsha Co. Ltd., Tokyo (1966) 170.

\bibitem{chen-xin} Q. Chen, Y. L. Xin, {\em A generalized maximum
principle and its applications in geometry.} Amer. J. of Math., {\bf
114}, 355--366, (1992).

\bibitem{cheng-yau}S. Y. Cheng and S. T. Yau, {\em Differential equations on Riemannian manifolds and their applications.} Comm. Pure. Appl. Math. \textbf{28} (1975) 333--354.

\bibitem{chern} S. S. Chern, {\em The Geometry of G-structures.} Bull. Amer. Math. Soc. {\bf 72} (1966), 167--219.

\bibitem{chern-kuiper} S. S. Chern and N. H. Kuiper, {\em Some theorems on the isometric imbedding of compact Riemannian manifolds in Euclidean space.} Ann. of Math. (2) \textbf{56} (1952) 442--430.

\bibitem{dias} C. C. Dias, {\em Isometric immersions with slow growth of curvature.} An. Acad. Bras. Ci\^{e}nc. {\bf 54}, 293--295, (1982).

\bibitem{efimov} N. Efimov,  {\em    Hyperbolic problems in the theory of surfaces.} Proc. Int. Congress Math. Moscow (1966);  Am. Math. Soc. translation \textbf{70} (1968), 26–38.

\bibitem{gimenez} F. Gim\'enez, {\em Estimates for the curvature of a submanifold which lies inside a tube}. J. Geom. {\bf 58} (1997), 95--105.

\bibitem{hasanis-koutroufiotis}  Th. Hasanis and Koutroufiotis, {\em Immersions of Riemannian manifolds into cylinders}. Arch. Math. {\bf 40} (1983), 82--85.

\bibitem{hilbert} D. Hilbert, {\em \"Uber Fl\"achen von konstanter Kr\"ummung.} Trans. Amer. Math. Soc. \textbf{2} (1901), 87-99.

\bibitem{jorge-koutrofiotis}L.  Jorge  and  D. Koutrofiotis,
{\em An estimate for the curvature of bounded submanifolds.} Amer. J. Math.  {\bf 103} (1980) 711--725.

\bibitem{jorge-xavier} L. Jorge, and F. Xavier, {\em An inequality between the exterior diameter and the mean
curvature of bounded immersions.} Math. Z. {\bf 178}, 77--82, (1981)

\bibitem{kobayashi-nomizu} S. Kobayashi, K. Nomizu, {\em Foundations of Differential Geometry,} vol II. Interscience Tracts in Pure and Appl. Math., no. 15, New York, (1969).

\bibitem{moore} J. D. Moore, {\em An application of second variation to submanifold theory.} Duke Math. J. \textbf{42}, (1975), 191--193.

\bibitem{nash} J. Nash, {\em The imbedding problem for Riemannian manifolds. } Ann. of Math. (2) \textbf{63},(1956) 20--63.

\bibitem{omori} H. Omori, {\em Isometric immersions of Riemannian manifolds.} J. Math. Soc. Japan \textbf{19}, (1967), 205--214.

\bibitem{oneil} B. O'Neill, {\em Immersions of manifolds of non-positive curvature.}  Proc. Amer. Math. Soc. \textbf{11} (1960), 132--134.

\bibitem{otsuki} T. Otsuki, {\em Isometric imbedding of Riemannian manifolds in a Riemannian manifold.} J. Math. Soc. Japan, \textbf{6}, (1954), 221--234.

\bibitem{pigola-rigoli-setti} S. Pigola, M. Rigoli and A. Setti,
{\em Maximum Principle on Riemannian Manifolds ans Applications.}
Mem. Amer. Math. Soc. {\bf 174}, no. 822 (2005).

\bibitem{schoen-yau} R. Schoen, \and S. T. Yau,  {\em Lectures on Differential Geometry.}
Conference Proceedings and Lecture Notes in Geometry and Topology,
 vol. \textbf{ 1}, (1994).

\bibitem{stiel} E. Stiel, {\em Immersions into manifolds of constant negative curvature.} Proc. Amer. Math. Soc. \textbf{18} (1967), 713--715.

\bibitem{tompkins}C. Tompkins, {\em Isometric embedding of flat manifolds in Euclidean spaces.} Duke Math. J. \textbf{5}, (1939), 58--61.

\bibitem{yau} S. T. Yau, {\em Harmonic functions on complete Riemannian manifolds.} Comm. Pure Appl. Math. \textbf{28}, (1975), 201--228.


\end{thebibliography}
\end{document}